# NONINTERSECTING BROWNIAN EXCURSIONS

By Craig A. Tracy[1] and Harold Widom[2]

*University of California-Davis and University of California-Santa Cruz*

We consider the process of $n$ Brownian excursions conditioned to be nonintersecting. We show the distribution functions for the top curve and the bottom curve are equal to Fredholm determinants whose kernel we give explicitly. In the simplest case, these determinants are expressible in terms of Painlevé V functions. We prove that as $n \to \infty$, the distributional limit of the bottom curve is the Bessel process with parameter $1/2$. (This is the Bessel process associated with Dyson's Brownian motion.) We apply these results to study the expected area under the bottom and top curves.

**1. Introduction.** A *Brownian excursion* (BE) $(X(\tau))_{\tau \in [0,1]}$ is a Brownian path conditioned to be zero at times $\tau = 0$ and $\tau = 1$ and to remain positive for $0 < \tau < 1$. This heuristic definition of this Markov process requires elaboration [11, 21] since such paths have zero Wiener measure. Alternatively, BE can be defined by scaling one-dimensional simple random walk conditioned to stay positive and conditioned to start and to end at the origin—a process known as *Bernoulli excursion* [22].

In this paper, we consider the process of $n$ Brownian excursions conditioned to be nonintersecting, a model first introduced by Gillet [14] and Katori et al. [16]. (Actually, the latter authors consider the closely related model of meanders instead of excursions.) Gillet proved that this $n$-nonintersecting BE is the distributional limit of suitably rescaled $n$ Bernoulli excursions conditioned not to intersect. This Bernoulli exclusion model is known in the physics literature as *vicious walkers with a wall* and was first proposed by Fisher [9] as a statistical physics model of wetting and melting. Using the theory of symmetric functions, Krattenthaler, Guttmann and Viennot

Received July 2006; revised November 2006.
[1]Supported by NSF Grant DMS-05-53379.
[2]Supported by NSF Grant DMS-05-52388.

[17] give explicit enumeration formulae for possible configurations of these walkers.

If the BE paths are labeled $X_1(\tau) < \cdots < X_n(\tau)$, $0 < \tau < 1$, then our main focus is on the lowest path, $X_1(\tau)$ and the highest path, $X_n(\tau)$. In particular, we show for each fixed positive integer $m$, with $0 < \tau_1 < \cdots < \tau_m < 1$, that

$$\mathbb{P}(X_1(\tau_1) \geq x_1, \ldots, X_1(\tau_m) \geq x_m)$$

and

$$\mathbb{P}(X_n(\tau_1) < x_1, \ldots, X_n(\tau_m) < x_m), \qquad x_i \in \mathbb{R}^+,$$

equal Fredholm determinants of the form $\det(I - K^{\mathrm{BE}}\chi_J)$, where $K^{\mathrm{BE}}$, the *extended BE kernel*, is an $m \times m$ matrix kernel and $\chi_J$ denotes multiplication by a diagonal matrix $\mathrm{diag}(\chi_{J_k})$ whose entries are indicator functions. For the lowest path, $J_k = (0, x_k)$ and for the highest path, $J_k = (x_k, \infty)$.

The matrix kernel $K^{\mathrm{BE}}(x, y) = (K^{\mathrm{BE}}_{k\ell}(x, y))_{k,\ell=1}^m$ will have entries given by differences,

$$K^{\mathrm{BE}}_{k\ell}(x, y) = H_{k\ell}(x, y) - E_{k\ell}(x, y).$$

Here,

(1.1) $$E_{k\ell}(x, y) = \begin{cases} P_-(x, y, \tau_\ell - \tau_k), & \text{if } k < \ell, \\ 0, & \text{otherwise}, \end{cases}$$

where $P_-(x, y, \tau)$ is the transition probability for Brownian excursion,

$$P_-(x, y, \tau) = \sqrt{\frac{2}{\pi\tau}} e^{-(x^2+y^2)/(2\tau)} \sinh\left(\frac{xy}{\tau}\right).$$

The summand $H_{k\ell}(x, y)$ will have different representations, one of which is

$$H_{k\ell}(x, y) = \sqrt{\frac{2}{\tau_\ell(1-\tau_k)}} e^{-x^2/(2(1-\tau_k)) - y^2/(2\tau_\ell)}$$

$$\times \sum_{j=0}^{n-1} \left[\frac{\tau_k(1-\tau_\ell)}{\tau_\ell(1-\tau_k)}\right]^{j+1/2}$$

$$\times p_{2j+1}\left(\frac{x}{\sqrt{2\tau_k(1-\tau_k)}}\right) p_{2j+1}\left(\frac{y}{\sqrt{2\tau_\ell(1-\tau_\ell)}}\right),$$

where the $p_k$ are the normalized Hermite polynomials. Although rather complicated, it has the same ingredients as the kernel for the Gaussian unitary ensemble (GUE).

Extended kernels such as these were first introduced by Eynard and Mehta [8] in the context of a certain Hermitian matrix model. In Section 2, we give a self-contained account of a class of extended kernels. Though the results appear in one form or another in the literature [8, 13, 26], we hope the reader will find this presentation beneficial. In Section 3, we use this theory of extended kernels to derive our expressions for the entries of $K^{\mathrm{BE}}(x,y)$.

Using the (scalar) $m = 1$ kernel, we show in Section 4 that for $\tau \in (0,1)$,

$$\mathbb{P}(X_1(\tau) \geq s\sqrt{2\tau(1-\tau)}) = \exp\left(-\int_0^{s^2} \frac{r(t)}{t}\, dt\right),$$

where $r = r(s)$ is the solution to

(1.2) $\quad s^2(r'')^2 = 4(r')^2(sr' - r) + (sr' - r)^2 - (4n+1)(sr' - r)r' + \tfrac{1}{4}(r')^2$

satisfying the asymptotic condition

$$r(s) = r_0 s^{3/2} + \mathrm{O}(s^{5/2}), \qquad s \to 0,$$

where

$$r_0 = \frac{1}{\sqrt{\pi}} \frac{1}{2^{2n}} \binom{2n}{n} \frac{4n(2n+1)}{3}.$$

A similar result is obtained for $\mathbb{P}(X_n(\tau) \leq s\sqrt{2\tau(1-\tau)})$.

The differential equations are the so-called $\sigma$ form of Painlevé V equations and the Fredholm determinants in these cases, where $J = (0,x)$ or $J = (x,\infty)$, are Painlevé V $\tau$-functions [12].

To obtain (1.2), we use results from [25] to obtain a system of differential equations whose solutions determine $r(s)$ and then reduce these to the second order equation for $r(s)$.

In Section 5, we consider the limiting process, as $n \to \infty$, of a suitably rescaled bottom path $X_1(\tau)$ near a fixed time $\tau \in (0,1)$ and show that the distributional limit is $\mathcal{B}_{1/2}(\tau)$, the *Bessel process* of [26] with parameter $\alpha = 1/2$. (This is different from the three-dimensional Bessel process appearing in the literature, which is Brownian motion conditioned to stay positive forever. When one rescales BE near $\tau = 0$, the three-dimensional Bessel process is the scaling limit.) This is a process whose finite-dimensional distribution functions are given in terms of an *extended Bessel kernel* $K^{\mathrm{Bes}}$. Precisely, if $-\infty < \tau_1 < \cdots < \tau_m < \infty$, then

(1.3) $\qquad \mathbb{P}(\mathcal{B}_{1/2}(\tau_1) \geq x_1, \ldots, \mathcal{B}_{1/2}(\tau_m) \geq x_m) = \det(I - K^{\mathrm{Bes}}\chi_J),$

where $\chi_J = \mathrm{diag}(\chi_{J_k})$ with $J_k = (0, x_k)$. For $\alpha = 1/2$, the entries of the matrix kernel are

$$K^{\mathrm{Bes}}_{k\ell}(x,y) = \begin{cases} \dfrac{2}{\pi} \displaystyle\int_0^1 e^{(\tau_k - \tau_\ell)t^2/2} \sin xt \sin yt\, dt, & \text{if } k \geq \ell, \\ -\dfrac{2}{\pi} \displaystyle\int_1^\infty e^{(\tau_k - \tau_\ell)t^2/2} \sin xt \sin yt\, dt, & \text{if } k < \ell. \end{cases}$$

Our result is that for fixed $\tau \in (0,1)$, if, in the probability

$$\mathbb{P}(X_1(\tau_1) \geq x_1, \ldots, X_1(\tau_m) \geq x_m),$$

we make the substitutions

$$x_k \to \sqrt{\frac{\tau(1-\tau)}{2n}} x_k, \qquad \tau_k \to \tau + \frac{\tau(1-\tau)}{2n} \tau_k,$$

then the resulting probability has the $n \to \infty$ limit (1.3). This is proved by establishing trace norm convergence of the extended kernels. As far as we know, this is the first appearance of the Bessel process in the literature as the distributional limit of a finite-$n$ Brownian motion process.

We expect the top path, $X_n(\tau)$, again suitably rescaled, to have the *Airy process* [13, 20] as its distributional limit. The reason is that the top path, as $n \to \infty$, should not feel the presence of the "wall" at $x = 0$ and the Airy process is known to be the distributional limit of nonintersecting Brownian bridges. The bottom path does feel the presence of the wall and in the end this is why the Bessel process appears. We shall not derive the scaling limit at the top because it is the bottom scaling that gives something new, and is thus more interesting.

For BE, a related random variable is the area $A$ under the BE curve. (See [10, 18] and references therein.) Both the distribution of $A$ and the moments of $A$ are known (the moments are known recursively) and, perhaps surprisingly, numerous applications have been found [10, 19] outside their original setting. It is natural, therefore, to introduce $A_{n,L}$ ($A_{n,H}$) equal to the area under the lowest (highest) curve in $n$-nonintersecting BE. From our limit theorem specialized to the $m = 1$ case, we deduce in Section 5.4 the asymptotics of the first moments $\mathbb{E}(A_{n,L})$ and $\mathbb{E}(A_{n,H})$. In particular,

$$\mathbb{E}(A_{n,L}) \sim \frac{c_L}{\sqrt{n}}, \qquad n \to \infty,$$

where

$$c_L = \frac{\pi\sqrt{2}}{16} \int_0^\infty \det(I - K^{\mathrm{Bes}} \chi_{(0,s)}) \, ds \simeq 0.682808.$$

Here, $K^{\mathrm{Bes}}$ is the $m = 1$ Bessel kernel,

$$K^{\mathrm{Bes}}(x,y) = \frac{2}{\pi} \int_0^1 \sin xt \sin yt \, dt = \frac{1}{\pi}\left(\frac{\sin(x-y)}{x-y} - \frac{\sin(x+y)}{x+y}\right).$$

## 2. Extended kernels.

2.1. *Initial and final points all distinct.* Although only a special case of the following derivation will be used in our discussion of nonintersecting BE, it might be useful for the future, and is no more difficult, to consider a rather general setting.

Suppose that $\{X_i(\tau)\}_{i=1}^n$ is a family of stationary Markov processes with continuous paths, with common transition probability density $P(x,y,\tau)$, conditioned to be nonintersecting, to begin at $a_1,\ldots,a_n$ at time $\tau = 0$ and to end at $b_1,\ldots,b_n$ at time $\tau = 1$. An *extended kernel* is a matrix kernel $K(x,y) = (K_{k\ell}(x,y))_{k\ell=1,\ldots,m}$ depending on $0 < \tau_1 < \cdots < \tau_m < 1$ with the following property. Given functions $f_k$, the expected value of

$$\prod_{k=1}^{m}\prod_{i=1}^{n}(1+f_k(X_i(\tau_k)))$$

is equal to $\det(I+Kf)$, where $f$ denotes multiplication by $\mathrm{diag}(f_k)$. In the special case where $f_k = -\chi_{J_k}$, this is the probability that for $k=1,\ldots,m$, no path passes through the set $J_k$ at time $\tau_k$. For BE, this gives the statement made in the third paragraph of the Introduction, once we have computed the extended kernel for it.

Here is how an extended kernel is obtained. Define the $n \times n$ matrix $A$ by

$$A_{ij} = P(a_i, b_j, 1),$$

define the $m \times m$ matrix kernel $E(x,y)$ by

$$(2.1) \qquad E_{k\ell}(x,y) = \begin{cases} P(x,y,\tau_\ell - \tau_k), & \text{if } k < \ell, \\ 0, & \text{otherwise,} \end{cases}$$

and if $A$ is invertible, define the $m \times m$ matrix kernel $H(x,y)$ by

$$(2.2) \qquad H_{k\ell}(x,y) = \sum_{i,j=1}^{n} P(x,b_j,1-\tau_k)(A^{-1})_{ji}P(y,a_i,\tau_\ell).$$

We shall show that $H - E$ is an extended kernel in the above sense.

From Karlin and McGregor [15] and the Markov property, we know the following. For paths starting at $x_{0,i} = a_i$ at time $\tau = 0$, the probability density that at times $\tau_k$ ($k = 1,\ldots,m+1$) the paths are at points $x_{ki}$ ($i=1,\ldots,n$) is a constant times

$$\prod_{k=0}^{m} \det(P(x_{k,i}, x_{k+1,j}, \tau_{k+1} - \tau_k))_{1 \le i,j \le n}.$$

In the above, we set $x_{m+1,j} = b_j, \tau_{m+1} = 1$. The expected value in question is obtained by multiplying

$$\prod_{k=1}^{m}\prod_{i=1}^{n}(1+f_k(x_{ki}))$$

by the probability density and integrating over all of the $x_{ki}$. We apply the general identity

$$\int \cdots \int \det(\varphi_j(x_k))_{j,k=1}^n \cdot \det(\psi_j(x_k))_{j,k=1}^n \, d\mu(x_1) \cdots d\mu(x_n)$$
$$= n! \det\left(\int \varphi_j(x)\psi_k(x) \, d\mu(x)\right)_{j,k=1}^n$$

successively to the variables $x_{1i}, \ldots, x_{mi}$ and find that the integral equals a constant times the determinant of the matrix with $i, j$ entry

$$\int \cdots \int P(a_i, x_1, \tau_1) \prod_{k=1}^{m-1} P(x_k, x_{k+1}, \tau_{k+1} - \tau_k) P(x_m, b_j, 1 - \tau_m)$$
$$\times \prod_{k=1}^m (1 + f_k(x_k)) \, dx_1 \cdots dx_m.$$

If we set $f = 0$, the matrix is $A$ by the semigroup property of $P(x, y, \tau)$ and the expected value is 1, so the normalizing constant is $(\det A)^{-1}$. Hence, if we replace $P(a_i, x, \tau)$ by

$$Q(a_i, x, \tau) = \sum_j (A^{-1})_{ij} P(a_j, x, \tau),$$

then the expected value equals the determinant of the matrix with $i, j$ entry

$$\int \cdots \int Q(a_i, x_1, \tau_1) \prod_{k=1}^{m-1} P(x_k, x_{k+1}, \tau_{k+1} - \tau_k) P(x_m, b_j, 1 - \tau_m)$$
$$\times \prod_{k=1}^m (1 + f_k(x_k)) \, dx_1 \cdots dx_m,$$

which equals the determinant of $I$ plus the matrix with $i, j$ entry

$$\int \cdots \int Q(a_i, x_1, \tau_1) \prod_{k=1}^{m-1} P(x_k, x_{k+1}, \tau_{k+1} - \tau_k) P(x_m, b_j, 1 - \tau_m)$$
$$\times \left[\prod_{k=1}^m (1 + f_k(x_k)) - 1\right] dx_1 \cdots dx_m.$$

The bracketed expression may be written as a sum of products,

$$\sum_{r \geq 1} \sum_{k_1 < \cdots < k_r} f_{k_1}(x_{k_1}) \cdots f_{k_r}(x_{k_r}).$$

Correspondingly, the integral is a sum of integrals. By integrating with respect to all $x_k$ with $k \neq k_1, \ldots, k_r$, we see that the corresponding integral is

equal to

$$\int \cdots \int Q(a_i, x_{k_1}, \tau_{k_1}) f_{k_1}(x_{k_1})$$
$$\times P(x_{k_1}, x_{k_2}, \tau_{k_2} - \tau_{k_1}) f_{k_2}(x_{k_2}) \cdots P(x_{k_r}, b_j, 1 - \tau_{k_r}) \, dx_{k_1} \cdots dx_{k_r}.$$

If we let $U_{k,\ell}$ be the operator with kernel $U_{k,\ell}(x,y) = P(x,y,\tau_l - \tau_k) f_\ell(y)$, then the above may be written as the single integral

$$\int Q(a_i, x, \tau_{k_1}) f_{k_1}(x) \cdot (U_{k_1,k_2} \cdots U_{k_{r-1},k_r} P(x, b_j, 1 - \tau_{k_r})) \, dx.$$

(If $r = 1$, we interpret the operator product to be the identity.) Replacing the index $k_1$ by $k$ and changing notation, we see that the sum of all of these equals

$$\int \sum_k Q(a_i, x, \tau_k) f_k(x) \left( \sum_{r \geq 0} \sum_{k_1, \ldots, k_r} U_{k,k_1} U_{k_1,k_2} \cdots U_{k_{r-1},k_r} P(x, b_j, 1 - \tau_{k_r}) \right) dx,$$

where the inner sum runs over all $k_r > \cdots > k_1 > k$. (Now, if $r = 0$, we interpret the operator product to be the identity.)

Write $P_j$ for the vector function $\{P(x, b_j, 1 - \tau_k)\}_{k=1}^m$ and $U$ for the operator acting on vector functions with matrix kernel $(U_{k,\ell}(x,y))_{k,\ell=1}^m$, where we set $U_{k,\ell}(x,y) = 0$ if $k \geq \ell$. Using the fact that $(I - U)^{-1} = \sum_{r \geq 0} U^r$, we see that the expression in large parentheses in the integral equals the $k$th component of the vector function $(I - U)^{-1} P_j$. Denote it by $T(k, x; j)$, and similarly denote the first factor in the integral by $S(i; k, x)$. We may think of $T(k, x; j)$ as the kernel of an operator $T$ taking sequences into vector functions,

$$\{a_j\}_{j=1}^n \to \left\{ \sum_{j=1}^n T(k, x; j) a_j \right\}_{k=1}^m,$$

while $S(i; k, x)$ is the kernel of an operator taking vector functions to sequences,

$$\{\varphi_k(x)\}_{k=1}^m \to \left\{ \int \sum_{k=1}^m S(i; k, x) \varphi_k(x) \, dx \right\}_{i=1}^n.$$

The integral is the kernel of the product $ST$, evaluated at $i, j$.

We use the general fact that $\det(I + ST) = \det(I + TS)$. The operator $TS$, which takes vector functions to vector functions, has matrix kernel with $k, \ell$ entry

$$\sum_{j=1}^n T(k, x; j) S(j; \ell, y).$$

Since $U$ is strictly upper-triangular, we have $\det(I - U) = 1$, so
$$\det(I + TS) = \det(I - U + (I - U)TS).$$
If we recall that $T(k, x; j)$ is the $k$th component of $(I - U)^{-1}$ applied to the vector function $\{P(x, b_j, 1 - \tau_k)\}$ and recall the definition of $S$, we see that $(I - U)TS$ has $k, \ell$ entry
$$\sum_{j=1}^{n} P(x, b_j, 1 - \tau_k) Q(a_j, y, \tau_\ell) f_\ell(y)$$
$$= \sum_{i,j=1}^{n} P(x, b_j, 1 - \tau_k)(A^{-1})_{ji} P(y, a_i, \tau_\ell) f_\ell(y).$$
This equals $H_{k\ell}(x, y) f_\ell(y)$. Since $U$ has $k, \ell$ entry $E_{k\ell}(x, y) f_\ell(y)$, we have shown that the expected value in question equals $\det(I + (H - E)f)$. Thus $H - E$ is an extended kernel.

2.2. *Initial and final points at zero.* In the preceding derivation, we assumed that $A$ was invertible. If any two $a_i$ are equal or any two $b_j$ are equal, this will certainly not hold. There are analogous, although more complicated, formulas if several of the $a_i$ or $b_j$ are equal. We consider here the simplest case where they are all equal to zero. We take (2.2) for the case when they are all different (when $A$ is expected to be invertible) and compute the limit as all of the $a_i$ and $b_j$ tend to zero.

There is a matrix function $D_0 = D_0(c_1, \ldots, c_n)$ such that for any smooth function $f$,
$$\lim_{c_i \to 0} D_0 \begin{pmatrix} f(c_1) \\ f(c_2) \\ \vdots \\ f(c_n) \end{pmatrix} = \begin{pmatrix} f(0) \\ f'(0) \\ \vdots \\ f^{(n-1)}(0) \end{pmatrix}.$$
Here, $\lim_{c_i \to c}$ is short for a certain sequence of limiting operations. The matrix $D_0'$ (the transpose of $D_0$) acts in a similar way on the right on row vectors.

Writing $D_a$ for $D_0(a_1, \ldots, a_n)$ and $D_b$ for $D_0(b_1, \ldots, b_n)$, we use
$$A^{-1} = D_b'(D_a A D_b')^{-1} D_a$$
and then take limits as all $a_i, b_j \to 0$. The matrix $D_a A D_b'$ becomes the matrix $B$ with $i, j$ entry
$$B_{ij} = \partial_{a_i}^{i-1} \partial_{b_j}^{j-1} P(a_i, b_j, 1)|_{a_i = b_j = 0},$$
the row vector $\{P(x, b_j, 1 - \tau_k)\}$ becomes the vector with $j$th component
$$P_j^b(x, 1 - \tau_k) = \partial_{b_j}^{j-1} P(x, b_j, 1 - \tau_k)|_{b_j = 0}$$

and the column vector $\{P(a_i, y, \tau_\ell)\}$ becomes the vector with $i$th component
$$P_i^a(y, \tau_l) = \partial_{a_i}^{i-1} P(a_j, y, \tau_\ell)|_{a_i=0}.$$
Then if the matrix $B$ is invertible, we obtain for the $H$-summand of the extended kernel, in this case,

$$(2.3) \qquad H_{k\ell}(x,y) = \sum_{i,j=1}^{n} P_j^b(x, 1-\tau_k)(B^{-1})_{ji} P_i^a(y, \tau_\ell).$$

One would like $B$ to be as simple as possible. Something one could have done before passing to the limit would be to take two sequences $\{\alpha_i\}$ and $\{\beta_j\}$ of nonzero numbers and in (2.2), to replace $A$ by the matrix $(\alpha_i A_{ij} \beta_j)$, replace the function $P(x, b_j, 1-\tau_k)$ by $P(x, b_j, 1-\tau_k)\beta_j$ and replace the function $P(a_i, y, \tau_\ell)$ by $\alpha_i P(a_i, y, \tau_\ell)$. The kernel does not change.

2.2.1. *An example.* We apply (2.3) to the case of $n$ nonintersecting Brownian bridges. The transition probability for Brownian motion is
$$P(x, y, \tau) = \frac{1}{\sqrt{2\pi\tau}} e^{-(x-y)^2/2\tau}$$
and we choose $\alpha_i = e^{a_i^2/2}$ and $\beta_j = e^{b_j^2/2}$. Then the modified $A$ equals $(e^{a_i b_j}/\sqrt{2\pi})$, the modified $B$ equals $\operatorname{diag}((i-1)!)/\sqrt{2\pi}$, the modified $B^{-1}$ equals $\sqrt{2\pi} \times \operatorname{diag}(1/(i-1)!)$ and

$$P_j^b(x, 1-\tau_k) = \frac{1}{\sqrt{2\pi(1-\tau_k)}} e^{-x^2/(2(1-\tau_k))}$$
$$\times \left[\frac{\tau_k}{2(1-\tau_k)}\right]^{(j-1)/2} H_{j-1}\left(\frac{x}{\sqrt{2\tau_k(1-\tau_k)}}\right),$$

$$P_i^a(y, \tau_\ell) = \frac{1}{\sqrt{2\pi\tau_\ell}} e^{-y^2/(2\tau_\ell)} \left[\frac{1-\tau_\ell}{2\tau_\ell}\right]^{(i-1)/2} H_{i-1}\left(\frac{y}{\sqrt{2\tau_\ell(1-\tau_\ell)}}\right),$$

where $H_k(x)$ are the Hermite polynomials satisfying
$$e^{-t^2+2tx} = \sum_{k\geq 0} \frac{t^k}{k!} H_k(x).$$

In terms of the harmonic oscillator wave functions,

$$(2.4) \qquad \varphi_k(x) = e^{-x^2/2} p_k(x) = \frac{1}{\sqrt{\pi^{1/2} 2^k k!}} e^{-x^2/2} H_k(x),$$

we have
$$H_{k\ell}(x,y) = \frac{1}{\sqrt{2(1-\tau_k)\tau_\ell}} e^{(1/2-\tau_k)X_k^2}$$
$$\times \sum_{j=0}^{n-1} \left[\frac{\tau_k(1-\tau_\ell)}{\tau_\ell(1-\tau_k)}\right]^{j/2} \varphi_j(X_k) \varphi_j(Y_\ell) e^{-(1/2-\tau_\ell)Y_\ell^2},$$

where, for notational simplicity, we write $X_k = x/\sqrt{2\tau_k(1-\tau_k)}$ and $Y_\ell = y/\sqrt{2\tau_\ell(1-\tau_\ell)}$. The matrix kernel $H(x,y) - E(x,y)$ is called the *extended Hermite kernel*. For $m=1$, $H_{11}(x,y)$ (with $\tau_1 = \tau$) reduces to the Hermite kernel $K_n^{\mathrm{GUE}}(X,Y)$. [We can neglect the two exponential factors since they do not change the value of the determinant. The factor $1/\sqrt{2\tau(1-\tau)}$ also disappears when we multiply by the differential $dY$.] [In the usual definition of the extended Hermite kernel, the times $\tau_k \in (0,1)$ are replaced by times $\hat{\tau}_k \in (0,\infty)$ defined by $\tau_k/(1-\tau_k) = e^{2\hat{\tau}_k}$ so that $H - E$ has $k,\ell$ entry $K_n^{\mathrm{GUE}}(X,Y;\hat{\tau}_k - \hat{\tau}_\ell)$, where

$$(2.5) \quad K_n^{\mathrm{GUE}}(x,y;\hat{\tau}) = \begin{cases} \displaystyle\sum_{j=0}^{n-1} e^{j\hat{\tau}} \varphi_j(x)\varphi_j(y), & \text{if } \hat{\tau} \geq 0, \\ -\displaystyle\sum_{j=n}^{\infty} e^{j\hat{\tau}} \varphi_j(x)\varphi_j(y), & \text{if } \hat{\tau} < 0. \end{cases}$$

The representation for $\hat{\tau} < 0$ follows from the Mehler formula ([7], Section 10.13 (22)).]

2.3. *General case.* In the recent literature, one finds investigations where not all starting and ending points are the same. In [2, 27], there is one starting point and two ending points and in [6], there are two starting points and two ending points.

The previous discussion can be extended to the case where there are $R$ different $a_i$ and $S$ different $b_j$. The different $a_i$, we label $a_1, \ldots, a_R$ and assume that $a_r$ occurs $m_r$ times. Similarly, the different $b_j$, we label $b_1, \ldots, b_S$ and assume that $b_s$ occurs $n_s$ times. Of course, $\sum m_r = \sum n_s = n$. It is convenient to use double indices $ir$ and $js$ where $i = 0, \ldots, m_r - 1$ and $j = 0, \ldots, n_s - 1$. We may think of $A$ as consisting of $m_r \times n_s$ blocks $A_{rs}$.

Once again, there is a matrix function $D_c = D_c(c_1, \ldots, c_n)$ such that

$$\lim_{c_i \to c} D_c \begin{pmatrix} f(c_1) \\ f(c_2) \\ \vdots \\ f(c_n) \end{pmatrix} = \begin{pmatrix} f(c) \\ f'(c) \\ \vdots \\ f^{(n-1)}(c) \end{pmatrix}$$

and similarly for the action of $D'_c$ on row vectors. We then write

$$A^{-1} = D'_b(D_a A D'_b)^{-1} D_a,$$

where $D_a$ and $D_b$ are now block diagonal matrices,

$$D_a = \mathrm{diag}(D_{a_r}), \qquad D_b = \mathrm{diag}(D_{b_s}),$$

substitute into (2.2) and take the appropriate limits. The matrix entry $A_{ir,js}$ is replaced by

$$B_{ir,js} = \partial_{a_{ir}}^{i-1} \partial_{b_{js}}^{j-1} P(a_{ir}, b_{js}, 1)|_{a_{ir}=a_r, b_{js}=b_s},$$

while the row vector $\{P(x, b_{js}, 1 - \tau_k)\}$ becomes the vector with $js$ component

$$P_{js}^{b_j}(x, 1 - \tau_k) = \partial_{b_{js}}^{j-1} P(x, b_{js}, 1 - \tau_k)|_{b_{js}=b_s}$$

and the column vector $\{P(a_{ir}, y, \tau_\ell)\}$ becomes the vector with $ir$ component

$$P_{ir}^{a_r}(y, \tau_\ell) = \partial_{a_{ir}}^{j-1} P(a_{ir}, y, \tau_\ell)|_{a_{ir}=a_r}.$$

Thus (2.3) is now replaced by

$$H_{k\ell}(x, y) = \sum_{ir, js} P_{js}^{b_j}(x, 1 - \tau_k)(B^{-1})_{js, ir} P_{ir}^{a_r}(y, \tau_\ell).$$

## 3. Nonintersecting Brownian excursions.

3.1. *Preliminary result.* Let $B_\tau$ denote standard Brownian motion (BM), so that

$$\mathbb{P}_x(B_\tau \in dy) = P(x, y, \tau) \, dy = \frac{1}{\sqrt{2\pi\tau}} e^{-(x-y)^2/(2\tau)} \, dy, \qquad \tau > 0, \ x, y \in \mathbb{R},$$

and the subscript $x$ on $\mathbb{P}$ indicates that the BM was started at $x$. Let

$$H_a = \inf\{s : s > 0, B_s = a\} = \text{hitting time of } a.$$

The basic result we use is [11, 21]

$$(3.1) \quad \mathbb{P}_x(B_\tau \in dy, H_0 > \tau) = (P(x, y, \tau) - P(x, -y, \tau)) \, dy =: P_-(x, y, \tau) \, dy.$$

Equation (3.1) gives the probability that a Brownian path starting at $x$ is in the neighborhood $dy$ of $y$ at time $\tau$ while staying positive for all times $0 < s < \tau$. This process has continuous sample paths [11]; therefore, by Karlin–McGregor [15], we can construct an $n$ nonintersecting path version and apply the formalism of Section 2. It is convenient to write

$$(3.2) \qquad P_-(x, y, \tau) = \sqrt{\frac{2}{\pi\tau}} e^{-(x^2+y^2)/(2\tau)} \sinh\left(\frac{xy}{\tau}\right).$$

[The theory of Section 2 applies equally well to the $n$ nonintersecting path version of reflected (elastic) BM, $Y_\tau := |B_\tau|$. In this case,

$$\mathbb{P}_x(|B_\tau| \in dy) = (P(x, y, \tau) + P(x, -y, \tau)) \, dy$$
$$= \sqrt{\frac{2}{\pi\tau}} e^{-(x^2+y^2)/(2\tau)} \cosh\left(\frac{xy}{\tau}\right).$$

We leave it as an exercise to construct the extended kernel.]

3.2. *Extended kernel for Brownian excursion.* Referring to (2.2) and the technique described after (2.3), we find that for distinct $a_i$ and $b_i$,

$$H_{k\ell}(x,y) = \sqrt{\frac{\pi}{2}} \sum_{i,j=1}^n P_-(x,b_i,1-\tau_k) e^{b_i^2/2} (A^{-1})_{ij} P_-(y,a_j,\tau_\ell) e^{a_j^2/2},$$

where $A_{ij} = \sinh(a_i b_j)$. The even order derivatives of $\sinh x$ all vanish at $x=0$, so to take $a_i, b_j \to 0$, we have to modify the procedure described in Section 2.2. (Otherwise, the matrix $B$ so obtained would be noninvertible.) We now use the fact that there is a matrix function $D_{c,-}$ such that for any smooth odd function $f$,

$$\lim_{c_i \to 0} D_{c,-} \begin{pmatrix} f(c_1) \\ f(c_2) \\ \vdots \\ f(c_n) \end{pmatrix} = \begin{pmatrix} f'(0) \\ f'''(0) \\ \vdots \\ f^{(2n-1)}(0) \end{pmatrix}.$$

Thus the $B$ matrix is

$$B_{ij} = \partial_{a_i}^{2i-1} \partial_{b_j}^{2j-1} \sinh(a_i b_j)|_{a_i=b_j=0} = (2j-1)! \delta_{ij}.$$

It follows that the limiting kernel has $k, \ell$ entry

(3.3) $$H_{k\ell}(x,y) = \sqrt{\frac{\pi}{2}} \sum_{j=0}^{n-1} \frac{1}{(2j+1)!} P_j(x, 1-\tau_k) P_j(y, \tau_\ell),$$

where

$$P_j(x,\tau) = [\partial_y^{2j+1} P_-(x,y,\tau) e^{y^2/2}]|_{y=0}$$

$$= \frac{1}{\sqrt{2\pi\tau}} e^{-x^2/(2\tau)}$$

$$\times [\partial_y^{2j+1} e^{-((1-\tau)/(2\tau))y^2 + xy/\tau} - \partial_y^{2j+1} e^{-((1-\tau)/(2\tau))y^2 - xy/\tau}]|_{y=0}$$

$$= \sqrt{\frac{2}{\pi\tau}} \left(\frac{1-\tau}{2\tau}\right)^{j+1/2} e^{-x^2/(2\tau)} H_{2j+1}\left(\frac{x}{\sqrt{2\tau(1-\tau)}}\right).$$

Thus the $H$-summand of the extended kernel can be written as

(3.4) $$H_{k\ell}(x,y) = \sqrt{\frac{2}{\tau_\ell(1-\tau_k)}} e^{-x^2/(2(1-\tau_k)) - y^2/(2\tau_\ell)}$$

$$\times \sum_{j=0}^{n-1} \left[\frac{\tau_k(1-\tau_\ell)}{\tau_\ell(1-\tau_k)}\right]^{j+1/2} p_{2j+1}\left(\frac{x}{\sqrt{2\tau_k(1-\tau_k)}}\right)$$

$$\times p_{2j+1}\left(\frac{y}{\sqrt{2\tau_\ell(1-\tau_\ell)}}\right).$$

Recall that $p_k$ are the normalized Hermite polynomials. It follows from the Christoffel–Darboux formula (see, e.g., [7], Section 10.3 (10)) for Hermite polynomials that

$$\sum_{j=0}^{n-1} p_{2j+1}(x)p_{2j+1}(y)$$
$$= \frac{\sqrt{n}}{2}\left\{\frac{p_{2n}(x)p_{2n-1}(y) - p_{2n-1}(x)p_{2n}(y)}{x-y}\right.$$
$$\left. + \frac{p_{2n}(x)p_{2n-1}(y) + p_{2n-1}(x)p_{2n}(y)}{x+y}\right\}.$$

Setting $\tau_k = \tau$, $X = x/\sqrt{2\tau(1-\tau)}$ and $Y = y/\sqrt{2\tau(1-\tau)}$, we see that the diagonal elements are

$$(3.5) \quad \begin{aligned} H_{kk}(x,y) &= \sqrt{\frac{n}{2\tau(1-\tau)}} e^{-x^2/(2(1-\tau)) - y^2/(2\tau)} \\ &\quad \times \left\{\frac{p_{2n}(X)p_{2n-1}(Y) - p_{2n-1}(X)p_{2n}(Y)}{X-Y}\right. \\ &\quad \left. + \frac{p_{2n}(X)p_{2n-1}(Y) + p_{2n-1}(X)p_{2n}(Y)}{X+Y}\right\} \\ &= \frac{e^{(1/2-\tau)X^2} e^{-(1/2-\tau)Y^2}}{\sqrt{2\tau(1-\tau)}} \{K_{2n}^{\text{GUE}}(X,Y) - K_{2n}^{\text{GUE}}(X,-Y)\}. \end{aligned}$$

More generally (recall the notation of Section 2.2.1),

$$(3.6) \quad H_{k\ell}(x,y) - E_{k\ell}(x,y) = \frac{e^{(1/2-\tau_k)X_k^2 - (1/2-\tau_\ell)Y_\ell^2}}{\sqrt{2\tau_\ell(1-\tau_k)}}\{K_{2n}^{\text{GUE}}(X_k,Y_\ell;\hat{\tau}_k - \hat{\tau}_\ell)$$
$$- K_{2n}^{\text{GUE}}(X_k,-Y_\ell;\hat{\tau}_k - \hat{\tau}_\ell)\},$$

where $K_n^{\text{GUE}}(x,y;\hat{\tau})$ is defined as in (2.5). (To establish this identity for $k < \ell$, one uses the (slightly modified) Mehler formula,

$$\sum_{j=0}^{\infty} q^{2j+1}\varphi_{2j+1}(x)\varphi_{2j+1}(y)$$
$$= \frac{1}{\sqrt{\pi(1-q^2)}} e^{-(1+q^2)/(2(1-q^2))(x^2+y^2)} \sinh\left(\frac{2qxy}{1-q^2}\right),$$

with $x,y \to X_k, Y_\ell$ and $q = [\tau_k(1-\tau_\ell)/\tau_\ell(1-\tau_k)]^{1/2}$.)

3.2.1. *Density of paths.* The density of paths at time $\tau$ is given by

$$\rho_n(x,\tau) = H_{11}(x,x)$$

$$(3.7) \qquad = \frac{1}{\sqrt{2\tau(1-\tau)}}\{K_{2n}^{\text{GUE}}(X,X) - K_{2n}^{\text{GUE}}(X,-X)\}$$

$$= \frac{2}{\sqrt{2\tau(1-\tau)}} \sum_{j=0}^{n-1} \varphi_{2j+1}(X)^2$$

and satisfies the normalization

$$\int_0^\infty \rho_n(x,\tau)\,dx = n.$$

For $n=1$, $\rho_n$ reduces to the well-known density for Brownian excursion.

## 4. Distribution of lowest and highest paths at time $\tau$.

4.1. *Fredholm determinant representations.* Recall that we assume the nonintersecting BE paths to be labeled so that for all $0 < \tau < 1$,

$$X_n(\tau) > \cdots > X_1(\tau) > 0.$$

We are interested in the distribution of the lowest path, $X_1$ and the highest path, $X_n$. From Section 2, we see that

$$(4.1) \quad \mathbb{P}(X_1(\tau) \geq x) = \mathbb{P}(X_1(\tau) \geq x, \ldots, X_n(\tau) \geq x) = \det(I - K\chi_{J_1}),$$
$$J_1 = (0,s),$$

and

$$(4.2) \quad \mathbb{P}(X_n(\tau) < x) = \mathbb{P}(X_1(\tau) < x, \ldots, X_n(\tau) < x) = \det(I - K\chi_{J_2}),$$
$$J_2 = (s,\infty),$$

where, in both cases,

$$(4.3) \qquad K(x,y) = K_{2n}^{\text{GUE}}(x,y) - K_{2n}^{\text{GUE}}(x,-y)$$

and $s = x/\sqrt{2\tau(1-\tau)}$. [Recall that the relation between $K$ and $K^{\text{BE}}$ is

$$K(x,y) = \sqrt{2\tau(1-\tau)}K^{\text{BE}}(\sqrt{2\tau(1-\tau)}x, \sqrt{2\tau(1-\tau)}y),$$

aside from the exponential factor appearing outside the curly brackets in (3.6). For the determinants, we may ignore this factor.]

The operator $K$ is given by

$$K = \sum_{j=0}^{n-1} \Psi_j \otimes \Psi_j,$$

TABLE 1
*The expected areas of lowest and highest curves for $n = 1, \ldots, 9$*

| $n$ | $\mathbb{E}(A_{n,L})$ | $\mathbb{E}(A_{n,H})$ |
|---|---|---|
| 1 | $\sqrt{\frac{\pi}{8}} = 0.626657\ldots$ | $\sqrt{\frac{\pi}{8}} = 0.626657\ldots$ |
| 2 | $\frac{5}{8}(\sqrt{2}-1)\sqrt{\pi} = 0.458859\ldots$ | $\frac{5}{8}\sqrt{\pi} = 1.107783\ldots$ |
| 3 | $0.380211\ldots$ | $1.479479\ldots$ |
| 4 | $0.331955\ldots$ | $1.791454\ldots$ |
| 5 | $0.298441\ldots$ | $2.065097\ldots$ |
| 6 | $0.273407\ldots$ | $2.311582\ldots$ |
| 7 | $0.253784\ldots$ | $2.537567\ldots$ |
| 8 | $0.237865\ldots$ | $2.747386\ldots$ |
| 9 | $0.224612\ldots$ | $2.944040\ldots$ |

where $\Psi_j(x) = \sqrt{2}\varphi_{2j+1}(x)$ and where $u \otimes v$ for functions $u$ and $v$ denotes the operator $f \to u(v, f)$. It follows that

$$(4.4) \qquad \det(I - K\chi_J) = \det(\delta_{jk} - (\Psi_j, \Psi_k))_{j,k=0}^{n-1},$$

where the inner product $(\cdot, \cdot)$ is taken over the set $J$. Thus (4.4) gives finite determinant representations of $\mathbb{P}(X_1(\tau) \geq x)$ and $\mathbb{P}(X_n(\tau) < x)$.

As an application of (4.4), consider the areas under the lowest curve and the highest curve,

$$A_{n,L} := \int_0^1 X_1(\tau)\, d\tau, \qquad A_{n,H} := \int_0^1 X_n(\tau)\, d\tau.$$

For $n = 1$, $A_{1,L} = A_{1,H}$ and the distribution of this random variable is called the *Airy distribution* (not to be confused with the Airy process!), the moments of which are known recursively (see [10] and references therein). The expected area under the lowest curve is

$$\mathbb{E}(A_{n,L}) = \int_0^1 \mathbb{E}(X_1(\tau))\, d\tau = \int_0^1 \int_0^\infty \mathbb{P}(X_1(\tau) \geq x)\, dx\, d\tau$$

$$(4.5) \qquad = \int_0^1 \sqrt{2\tau(1-\tau)}\, d\tau \int_0^\infty \det(I - K\chi_{(0,s)})\, ds$$

$$(4.6) \qquad = \frac{\pi}{4\sqrt{2}} \int_0^\infty \det(\delta_{j,k} - (\Psi_j, \Psi_k))_{j,k=0}^{n-1}\, ds,$$

where the inner product is over the interval $(0, s)$. A similar formula holds for $\mathbb{E}(A_{n,H})$. Table 1 gives some data on $\mathbb{E}(A_{n,L})$ and $\mathbb{E}(A_{n,H})$. Higher moments are expressible in terms of integrals over the Fredholm determinants of the extended kernels. We have not investigated these higher moments.

4.2. *Integrable differential equations.* In this section, we show that the Fredholm determinant of the integral operator with kernel

(4.7) $$K(x,y)\chi_J(y) = [K_{2n}^{\text{GUE}}(x,y) - K_{2n}^{\text{GUE}}(x,-y)]\chi_J(y)$$

can be expressed in terms of a solution to a certain integrable differential equation when $J = (0,s)$ and $J = (s,\infty)$. To do this, we transform the kernel $K(x,y)$ to the canonical form

(4.8) $$\frac{\varphi(x)\psi(y) - \psi(x)\varphi(y)}{x-y}$$

so that the theory of Fredholm determinants and integrable differential equations, as developed in [25], can be directly applied. If, instead of $K(x,y)$, we use the kernel

$$K_0(x,y) = \frac{K(\sqrt{x},\sqrt{y})}{2(xy)^{1/4}},$$

then

(4.9) $$\det(I - K\chi_{(0,s^2)}) = \det(I - K_0\chi_{(0,s)}).$$

By use of the Christoffel–Darboux formula, we see that $K_0$ is of the form (4.8), with

(4.10) $$\varphi(x) = n^{1/4}x^{1/4}\varphi_{2n}(\sqrt{x}), \qquad \psi(x) = n^{1/4}x^{-1/4}\varphi_{2n-1}(\sqrt{x}),$$

where the $\varphi_k$ are the harmonic oscillator wave functions (2.4).

Using the well-known formulas for the action of multiplication by $x$ and differentiation by $x$ on the $\varphi_k$, we find

(4.11) $$x\varphi'(x) = (\tfrac{1}{4} - \tfrac{1}{2}x)\varphi(x) + \sqrt{n}x\psi(x),$$

(4.12) $$x\psi'(x) = -\sqrt{n}\varphi(x) - (\tfrac{1}{4} - \tfrac{1}{2}x)\psi(x).$$

These formulas will allow us to use the results of [25].

In the derivation of the next section, we write $J$ for $(0,s)$ and denote by $R(x,y)$ the resolvent kernel of $K_0\chi_J$, the kernel of $(I - K_0\chi_J)^{-1}K_0\chi_J$. Of course, this also depends on the parameter $s$. A basic fact is that

(4.13) $$\frac{d}{ds}\log\det(I - K_0\chi_J) = -R(s,s).$$

Once $R(s,s)$ is known, integration and exponentiation give the determinant.

4.2.1. *Reduction to Painlevé V, $J = (0, s)$.* Following [25], we define the functions of $x$ (which also depend on $s$ since $J$ does)

$$Q = (I - K_0 \chi_J)^{-1} \varphi, \qquad P = (I - K_0 \chi_J)^{-1} \psi$$

and then the quantities (which depend on $s$ only)

$$q = Q(s), \qquad p = P(s), \qquad r = sR(s,s),$$
$$u = (Q, \varphi), \qquad v = (Q, \psi), \qquad w(s) = (P, \psi),$$

where the inner products are taken over $J$.

The definitions are not important here. What is important are the relations among the various quantities, as established in [25]. Of course, we are interested in $r$. This function and its derivative are given in terms of the other functions by

(4.14)
$$r(s) = (w + \sqrt{n})q^2 + (\tfrac{1}{2} - s - 2\sqrt{n}w)qp$$
$$+ (-u + \sqrt{n}s + 2\sqrt{n}v)p^2,$$

(4.15)
$$\frac{dr}{ds} = -qp + \sqrt{n}p^2.$$

((4.14) is the specialization of (2.27) of [25] with additional insertions described in (2.32); (4.15) is the specialization of (2.28) of [25] with additional insertions described in (2.35). In both cases, we use the formulas (4.11) and (4.11).)

Specializing the differential equations of [25] to the case at hand, we find that

(4.16) $$s\frac{dq}{ds} = (\tfrac{1}{4} - \tfrac{1}{2}s - \sqrt{n}w)q + (-u + 2\sqrt{n}v + \sqrt{n}s)p,$$

(4.17) $$s\frac{dp}{ds} = -(w + \sqrt{n})q - (\tfrac{1}{4} - \tfrac{1}{2}s - \sqrt{n}w)p,$$

(4.18) $$\frac{du}{ds} = q^2, \qquad \frac{dv}{ds} = pq, \qquad \frac{dw}{ds} = p^2.$$

((4.16) and (4.17) are (2.25) and (2.26) of [25], resp., with additional insertions described in (2.31). We also used recursion relations (2.12) and (2.13) of [25], evaluated at $x = s$. Equations (4.18) are "universal equations" given by (2.16)–(2.18) of [25].)

This is a system of five equations in five unknowns. One has the first integrals, discovered by a computer search,

(4.19) $$sp(q - \sqrt{n}p) + \sqrt{n}(u + \tfrac{1}{2}w) - (1 + 2n)v + w(u - 2\sqrt{n}v + nw) = 0,$$

(4.20)
$$(\sqrt{n}s - u + 2\sqrt{n}v)p^2 + (\sqrt{n} + w)q^2$$
$$+ (\tfrac{1}{2} - s - 2\sqrt{n}w)pq + v - \sqrt{n}w = 0.$$

(A *first integral* for a system of differential equations is a relation stating that a particular expression in the variables is constant when the system is satisfied. Often, as here, the constant can be determined by using the initial conditions.)

To verify either of these, we differentiate the left-hand side using (4.16)–(4.18) and find that the derivative equals zero. Thus the left-hand side is a constant. That the constant is zero follows from the easily verified fact that all quantities are $o(1)$ as $s \to 0$.

In addition to (4.14) and (4.15), we compute $r''$ using (4.16)–(4.18). These, together with the first integrals (4.19) and (4.20), give us five equations. Remarkably, from these, we can eliminate $q$, $p$, $u$, $v$ and $w$ and we obtain the second order equation (1.2) for $r$. (For the elimination, we used a Gröbner basis algorithm as implemented in MATHEMATICA. A Gröbner basis of an ideal in a polynomial ring is a particular basis that is well suited for the elimination of variables. The reader is referred to the book by Cox, Little and O'Shea [5] for both the theory of Gröbner bases and many examples of their application to computational commutative algebra. Remarkably, although one would expect to need six equations to eliminate the five unknowns, we are in a nongeneric situation and these five equations suffice. This would not have been discovered without the use of the Gröbner basis. In our original derivation, we used only the first integral (4.19) and obtained a third order equation in $r$. Cosgrove [4] knew how to reduce this to the second order equation (1.2). This then led us to discover (4.20).) The solution we want satisfies the asymptotic condition

$$(4.21) \qquad r(s) = r_0 s^{3/2} + \mathrm{O}(s^{5/2}) \qquad \text{as } s \to 0,$$

where

$$(4.22) \qquad r_0 = \frac{1}{\sqrt{\pi}} \frac{1}{2^{2n}} \binom{2n}{n} \frac{4n(2n+1)}{3}.$$

This condition follows from the fact that $R(s,s) \sim K_0(s,s)$ as $s \to 0$. (Power counting suggests that $cs^{1/2}$ is the leading term, but a computation shows that $c$ is zero. The $s^{3/2}$ term is the leading one. Higher terms can be determined from the differential equation.)

4.2.2. *Painlevé V representation of* $\mathbb{P}(X_1(\tau) \geq x)$. From (4.13) and (4.9), we conclude that

$$(4.23) \qquad \mathbb{P}(X_1(\tau) \geq s\sqrt{2\tau(1-\tau)}) = \exp\left(-\int_0^{s^2} \frac{r(t)}{t}\,dt\right),$$

where $r$ is the solution to (1.2) satisfying boundary condition (4.21). Using (1.2), it is routine to compute the small $s$ expansion of $r(s)$ and hence that

of (4.23),

$$\mathbb{P}(X_1(\tau) \geq s\sqrt{2\tau(1-\tau)}) = 1 - \frac{2}{3}r_0 s^3 + \frac{2}{25}(4n+1)r_0 s^5 - \frac{64n^2+32n+9}{735}r_0 s^7$$
$$+ \frac{(4n+1)(32n^2+16n+15)}{8505}r_0 s^9$$
$$+ \frac{128(2n+3)(n-1)}{275625}r_0^2 s^{10} + \mathrm{O}(s^{11}),$$

where $r_0$ is defined as in (4.22).

We remark that the number of BE curves, $n$, appears only as a parameter in the differential equation (1.2) and the boundary condition (4.21). Although we have explicit solutions in terms of finite $n$ determinants, expansions such as the one above are more easily derived from the differential equation representation.

4.2.3. *Reduction to Painlevé V, $J = (s, \infty)$.* This differs from $J = (0, s)$ in only minor ways. In this case,

$$\frac{d}{ds} \log \det(I - K_0 \chi_J) = R(s,s),$$

so that with $r$ defined as before,

$$\mathbb{P}(X_n(\tau) \leq s\sqrt{2\tau(1-\tau)}) = \exp\left(-\int_{s^2}^{\infty} \frac{r(t)}{t}\,dt\right).$$

Specializing the differential equations of [25], we find for $J = (s, \infty)$,

$$s\frac{dq}{ds} = (\tfrac{1}{4} - \tfrac{1}{2}s - \sqrt{n}w)q + (-u + 2\sqrt{n}v + \sqrt{n}s)p,$$
$$s\frac{dp}{ds} = -(w + \sqrt{n})q - (\tfrac{1}{4} - \tfrac{1}{2}s - \sqrt{n}w)p,$$
$$\frac{du}{ds} = -q^2, \qquad \frac{dv}{ds} = -pq, \qquad \frac{dw}{ds} = -p^2.$$

The connection between $r$ and the above variables remains the same as that given in (4.14). One has the two first integrals

$$-sp(q - \sqrt{n}p) + \sqrt{n}(u + \tfrac{1}{2}w) - (1 + 2n)v + w(u - 2\sqrt{n}v + nw) = 0,$$
$$(\sqrt{n}s - u + 2\sqrt{n}v)p^2 + (\sqrt{n} + w)q^2 + (\tfrac{1}{2} - s - 2\sqrt{n}w)pq - v + \sqrt{n}w = 0.$$

Proceeding with the calculations as described above, we obtain

$$s^2(r'')^2 = -4(r')^2(sr' - r) + (sr' - r)^2 - (4n+1)(sr' - r)r' + \tfrac{1}{4}(r')^2.$$

**5. Limit theorems as $n \to \infty$.** In this section, we obtain the scaling limit for the extended BE kernel at the bottom. We first explain what this means and then state the result.

Given a kernel $L(x,y)$, the kernel resulting from substitutions $x \to \varphi(x), y \to \varphi(y)$ is the *scaled kernel*

$$\hat{L}(x,y) = (\varphi'(x)\varphi'(y))^{1/2} L(\varphi(x), \varphi(y)),$$

which defines a unitarily equivalent operator. If there is an underlying parameter and $\hat{L}$ has a limit $\hat{L}_0$ (pointwise or in some norm), then $\hat{L}_0$ is the *scaling limit* of $L(x,y)$. If the convergence is in trace norm, then the limit of $\det(I-L)$ equals $\det(I-\hat{L}_0)$. Trace norm convergence is assured if the underlying interval is bounded and if $\hat{L}(x,y)$ and its first partial derivatives converge uniformly.

The extended Bessel kernel $K^{\mathrm{BES}}$ [24], Section 8, is the matrix kernel whose entries are given by

$$K_{k\ell}^{\mathrm{BES}}(x,y) = \begin{cases} \displaystyle\int_0^1 e^{z^2(\tau_k-\tau_\ell)/2} \Phi_\alpha(xz)\Phi_\alpha(yz)\,dz, & \text{if } k \geq \ell, \\ \displaystyle-\int_1^\infty e^{z^2(\tau_k-\tau_l)/2} \Phi_\alpha(xz)\Phi_\alpha(yz)\,dz, & \text{if } k < \ell, \end{cases}$$

where $\Phi_\alpha(z) = \sqrt{z}J_\alpha(z)$. In particular, $\Phi_{1/2}(z) = \sqrt{2/\pi}\sin z$.

We shall show the following. Take a fixed $\tau \in (0,1)$ and make the substitutions

$$(5.1) \qquad x \to \sqrt{\frac{\tau(1-\tau)}{2n}} x, \qquad y \to \sqrt{\frac{\tau(1-\tau)}{2n}} y$$

and replacements

$$(5.2) \qquad \tau_k \to \tau + \frac{\tau(1-\tau)}{2n}\tau_k.$$

Then the extended BE kernel has the extended Bessel kernel with parameter $\alpha = 1/2$ as scaling limit as $n \to \infty$, in trace norm over any bounded domain. In particular, we deduce the result stated in the introduction, that the finit-dimensional distribution functions for the BE process, when suitably normalized, converge to the corresponding distribution functions for the Bessel process.

5.1. *Integral representation for the extended kernel.* The extended BE kernel $K^{\mathrm{BE}}$ has the form $(H_{k\ell}) - (E_{k\ell})$, where $E_{k\ell}$ is given by (1.1) and $H_{k\ell}$ is given by (3.3). The limit theorem for the $H_{k\ell}$ will be derived from an integral representation of (3.3) which, as in [3, 27], uses two different integral representations for $P_j$. (We could use known integral representations for the

Hermite polynomials, but the method we use here is more general and also works in other cases.)

The first is easy. With integration over a curve around zero,
$$P_-(x,y,\tau)e^{y^2/2} = \frac{1}{2\pi i}\int P_-(x,t,\tau)e^{t^2/2}\frac{dt}{t-y},$$
so
$$\begin{aligned}P_j(x,\tau) &= \frac{(2j+1)!}{2\pi i}\oint P_-(x,t,\tau)e^{t^2/2}\frac{dt}{t^{2j+2}}\\ &= \sqrt{\frac{2}{\pi\tau}}e^{-x^2/2\tau}\frac{(2j+1)!}{2\pi i}\oint e^{-(1-\tau)/(2\tau)t^2}\sinh\frac{xt}{\tau}\frac{dt}{t^{2j+2}}.\end{aligned}$$

For the other, we take the Fourier transform on the $y$-variable. A computation gives
$$\begin{aligned}\int_{-\infty}^{\infty}e^{isy}P_-(x,y,\tau)&e^{y^2/2}\,dy\\ &= \frac{2}{\sqrt{1-\tau}}e^{x^2/2(1-\tau)}e^{-\tau/(2(1-\tau))s^2}\sinh\frac{ixs}{1-\tau}.\end{aligned}$$
Taking inverse Fourier transforms and differentiating gives
$$\begin{aligned}P_j(x,\tau) &= \frac{2}{\pi\sqrt{1-\tau}}e^{x^2/2(1-\tau)}\\ &\quad\times\int_{-\infty}^{\infty}e^{-\tau/(2(1-\tau))s^2}\sinh\frac{ixs}{1-\tau}(-is)^{2j+1}\,ds\\ &= \frac{2i}{\pi\sqrt{1-\tau}}e^{x^2/2(1-\tau)}\int_{-i\infty}^{i\infty}e^{\tau/(2(1-\tau))s^2}\sinh\frac{xs}{1-\tau}s^{2j+1}\,ds.\end{aligned}$$
These give
$$\begin{aligned}H_{k\ell}(x,y) &= \frac{1}{\pi^2\sqrt{(1-\tau_k)(1-\tau_\ell)}}e^{-x^2/(2(1-\tau_k))+y^2/(2(1-\tau_\ell))}\\ &\quad\times\oint\int_{-i\infty}^{i\infty}e^{-\tau_k/(2(1-\tau_k))t^2+\tau_\ell/(2(1-\tau_\ell))s^2}\sinh\frac{xt}{1-\tau_k}\\ &\quad\times\sinh\frac{ys}{1-\tau_\ell}\left[\left(\frac{s}{t}\right)^{2n}-1\right]\frac{s\,ds\,dt}{s^2-t^2}.\end{aligned}$$
If we ensure that the $s$-contour passes to one side of the $t$-contour, then the summand 1 in the bracket contributes zero since the $t$-integral is zero. Thus we have
$$H_{k\ell}(x,y) = \frac{1}{\pi^2\sqrt{(1-\tau_k)(1-\tau_\ell)}}e^{-x^2/(2(1-\tau_k))+y^2/(2(1-\tau_\ell))}$$

(5.3)
$$\times \oint \int_{-i\infty}^{i\infty} e^{-\tau_k/(2(1-\tau_k))t^2 + \tau_\ell/(2(1-\tau_\ell))s^2}$$

$$\times \sinh \frac{xt}{1-\tau_k} \sinh \frac{ys}{1-\tau_\ell} \left(\frac{s}{t}\right)^{2n} \frac{s\,ds\,dt}{s^2-t^2}.$$

5.2. *Limit theorem for the bottom.* If we ignore the denominator $s^2 - t^2$ in (5.3), we have a product of two integrals for each of which we can apply the method of steepest descent. If we think of $x$ and $y$ as small and $k = \ell$ with $\tau_k = \tau_\ell = \tau$, then there are saddle points at $\pm i\sqrt{2n(1-\tau)/\tau}$ for both integrals. So we assume that all $\tau_k$ are near a single $\tau$ and make the changes of variable

$$s \to \sqrt{\frac{2n(1-\tau)}{\tau}} s, \qquad t \to \sqrt{\frac{2n(1-\tau)}{\tau}} t.$$

This is what suggested the substitutions (5.1) and replacements (5.2). Observe that in the scaling limit, $E_{k,\ell}$ remains the same.

For (5.3), we can ignore the outside factor involving $x$ and $y$ since removing it does not change the determinant. An easy computation shows that after the various substitutions, the rest of $H_{k\ell}(x,y)$ becomes

$$\frac{1}{\pi^2 + O(n^{-1})} \oint \int_{-i\infty}^{i\infty} e^{n(s^2-t^2+2\log s - 2\log t) - (\tau_k + O(n^{-1}))t^2/2 + (\tau_\ell + O(n^{-1}))s^2/2}$$

$$\times \sinh \frac{xt}{1+O(n^{-1})} \sinh \frac{ys}{1+O(n^{-1})} \frac{s\,ds\,dt}{s^2 - t^2},$$

where the terms $O(n^{-1})$ are independent of the variables.

The saddle points are now at $\pm i$ for both integrals. The steepest descent curve for the $s$-integral would be the imaginary axis and for the $t$-integral would be the lines $\operatorname{Im} t = \pm 1$, the lower one described left to right, the upper one right to left. But in order to get to these steepest descent curves, we must pass the $s$-contour through part of the original $t$-contour. To be specific, assume that the $t$-contour was the unit circle described counterclockwise and that the $s$-contour passed to the right of the $t$-contour. Then if we replace the $s$-contour by the imaginary axis (after which we can replace the $t$-contour by the lines $\operatorname{Im} t = \pm 1$), we must add the contribution of the pole at $s = t$ when $t$ is in the right half-plane. Thus the change of contour results in having to add

$$\frac{-1}{\pi i + O(n^{-1})} \int_{-i}^{i} e^{(\tau_\ell - \tau_k + O(n^{-1}))t^2/2} \sinh \frac{xt}{1+O(n^{-1})} \sinh \frac{yt}{1+O(n^{-1})} dt.$$

This is

$$\frac{2}{\pi} \int_0^1 e^{(\tau_k - \tau_\ell)t^2/2} \sin xt \sin yt\,dt + O(n^{-1}),$$

uniformly for $x$ and $y$ bounded.

It is easy to see that the remaining integral, the $s$-integral over the imaginary axis and the $t$ integral over $\operatorname{Im} t = \pm 1$, is $O(n^{-1/2})$ uniformly for $x, y$ bounded. In fact, outside neighborhoods of the points where $s^2 = t^2$, we may estimate the double integral by a product of single integrals, each of which is $O(n^{-1/2})$, so the product is $O(n^{-1})$. There are four points where $s^2 = t^2$, where both $s$ and $t$ are $\pm i$ or where one is $i$ and the other $-i$. If we set up polar coordinates in neighborhoods based at any of these points, then the integrand, aside from the last factor, is $O(e^{-\delta n r^2})$ for some $\delta > 0$, while the last factor is $O(r^{-1})$ because $s$ and $t$ (or $s$ and $-t$) are on mutually perpendicular lines. Therefore, because of the factor $r$ coming from the element of area in polar coordinates, the integral is $O(\int_0^\infty e^{-\delta n r^2} dr) = O(n^{-1/2})$.

Thus $H_{k\ell}(x, y)$ has the scaling limit

$$\frac{2}{\pi} \int_0^1 e^{(\tau_k - \tau_\ell) t^2 / 2} \sin xt \sin yt \, dt$$

uniformly for $x$ and $y$ bounded. When $k < \ell$, we must subtract from this

$$E_{k\ell}(x, y) = P_-(x, y, \tau_\ell - \tau_k),$$

which equals the expression above, but with the integral taken over $(0, \infty)$. So the effect of the subtraction is to replace the expression by

$$-\frac{2}{\pi} \int_1^\infty e^{(\tau_k - \tau_\ell) t^2 / 2} \sin xt \sin yt \, dt.$$

This limiting matrix kernel is exactly the extended Bessel kernel $K^{\mathrm{Bes}}(x, y)$ with $\alpha = 1/2$.

We have shown that the $k, \ell$ entry of the scaled BE kernel converges to the corresponding entry of the extended Bessel kernel, uniformly for $x, y$ bounded. A minor modification of the computation shows that the corresponding result holds for derivatives with respect to $x$ and $y$ of all orders and this is more than sufficient to give trace norm convergence on any bounded set. This concludes the proof of the result stated at the beginning of the section.

5.3. *Limit theorem for highest curve.* As mentioned in the Introduction, we shall not derive the scaling limit at the top. But since below, we shall use the top scaling of the diagonal terms $H_{kk}(x, y)$ to the Airy kernel, we point out that this is very easy. With the usual scaling $X \to \sqrt{2n} + X/(2^{1/2} n^{1/6})$, $Y \to \sqrt{2n} + Y/(2^{1/2} n^{1/6})$, the summand in (3.5) with the denominator $X + Y$ is easily seen to converge in trace norm to zero over any interval $(s, \infty)$, by the asymptotics of the Hermite polynomials. The rest of the kernel converges to the Airy kernel in trace norm, by known results.

5.4. *Asymptotics of $\mathbb{E}(A_{n,L})$ and $\mathbb{E}(A_{n,H})$.* From our limit theorems follow the asymptotics of $\mathbb{E}(A_{n,L})$ and $\mathbb{E}(A_{n,H})$.

For the lowest curve, we use formula (4.5), which tells us that

$$\mathbb{E}(A_{n,L}) = \frac{\pi}{4\sqrt{2}} \int_0^\infty \det(I - K\chi_{(0,s)})\,ds,$$

where $K$ is given by (4.3). It is given in terms of $K^{\text{BE}}$ with $m=1$ and arbitrary $\tau_1 = \tau$ by

$$K(x,y) = \sqrt{2\tau(1-\tau)} K^{\text{BE}}(\sqrt{2\tau(1-\tau)}x, \sqrt{2\tau(1-\tau)}y)$$

times exponential factors which we may ignore. Hence, the scaling we have derived shows that

$$\frac{1}{2\sqrt{n}} K\left(\frac{x}{2\sqrt{n}}, \frac{y}{2\sqrt{n}}\right) \to K^{\text{Bes}}(x,y) = \frac{1}{\pi}\left(\frac{\sin(x-y)}{x-y} - \frac{\sin(x+y)}{x+y}\right)$$

in trace norm on bounded sets. As in Section 4.2, instead of $K^{\text{Bes}}$, we use $K_0^{\text{Bes}}$ defined by

$$K_0^{\text{Bes}}(x,y) = \frac{K^{\text{Bes}}(\sqrt{x}, \sqrt{y})}{2x^{1/4}y^{1/4}}.$$

It is convenient to use $K_0^{\text{Bes}}$ since we can make direct use of the differential equations for it derived in [24].

By a change of variable, we have

$$\mathbb{E}(A_{n,L}) = \frac{\pi}{8\sqrt{2n}} \int_0^\infty \det(I - K\chi_{(0,s/2\sqrt{n})})\,ds.$$

Therefore, if we can take the limit under the integral sign (which we shall show below), then this implies that

(5.4) $$\mathbb{E}(A_{n,L}) \sim \frac{c_L}{\sqrt{n}}, \qquad n \to \infty,$$

where

$$c_L = \frac{\pi}{8\sqrt{2}} \int_0^\infty \det(I - K_0^{\text{Bes}}\chi_{(0,s^2)})\,ds = \frac{\pi}{16\sqrt{2}} \int_0^\infty \frac{1}{\sqrt{x}} \det(I - K_0^{\text{Bes}}\chi_{(0,x)})\,dx.$$

After expressing this in terms of $K^{\text{Bes}}$, it becomes the formula quoted in the Introduction.

Solving numerically the differential equation satisfied by the logarithmic derivative of $\det(I - K_0^{\text{Bes}}\chi_{(0,x)})$ [24], followed by a numerical integration, gives

$$\int_0^\infty \frac{1}{\sqrt{x}} \det(I - K_0^{\text{Bes}}\chi_{(0,x)})\,dx \simeq 4.917948$$

NONINTERSECTING BROWNIAN EXCURSIONS 25and hence $c_L \simeq 0.682808$. From Table 1, we compute $\sqrt{n}\mathbb{E}(A_{n,L})$ for $n = 5,\ldots,9$:

$$0.667334, \quad 0.669708, \quad 0.671449, \quad 0.672784, \quad 0.673838.$$

To justify taking the limit under the integral sign in

$$\int_0^\infty \det(I - K\chi_{(0,s/2\sqrt{n})})\,ds,$$

it is enough to show two things:

(5.5) $$\lim_{n\to\infty} \int_{\sqrt{n}\log n}^\infty \det(I - K\chi_{(0,s/\sqrt{n})})\,ds = 0,$$

(5.6) $$\det(I - K\chi_{(0,s/2\sqrt{n})}) = O(e^{-\delta s}) \text{ uniformly for some } \delta \text{ when } s < \sqrt{n}\log n.$$

Combining (5.5) with an application of the dominated convergence theorem using (5.6) justifies taking the limit under the integral sign. This completes the proof of (5.4).

The eigenvalues of $K\chi_{(0,s/2\sqrt{n})}$ lie in the interval $[0,1]$ because $K$ is a projection operator. For (5.5), we use the fact that the largest eigenvalue of $K\chi_{(0,s/2\sqrt{n})}$ is at least that of $2\varphi_1 \otimes \varphi_1 \chi_{(0,s/\sqrt{n})}$, so

$$\det(I - K\chi_{(0,s/2\sqrt{n})}) \leq 1 - 2\int_0^{s/2\sqrt{n}} \varphi_1(x)^2\,dx = \int_{s/2\sqrt{n}}^\infty \varphi_1(x)^2\,dx = O(e^{-s^2/4n}).$$

Integrating this over $(\sqrt{n}\log n, \infty)$ gives (5.5).

For (5.6), we use the fact that since the spectrum of $K\chi_{(0,s/2\sqrt{n})}$ lies in $(0,1)$, we have

$$\det(I - K\chi_{(0,s/2\sqrt{n})}) \leq e^{-\mathrm{tr} K\chi_{(0,s/2\sqrt{n})}}.$$

(This follows from the inequality $1 - \lambda \leq e^{-\lambda}$ when $\lambda \in [0,1]$.) The trace equals

$$\int_0^{s/2\sqrt{n}} K(x,x)\,dx.$$

Now, the usual asymptotics for the Hermite polynomials for "small" $x$ hold uniformly for $x = O(n^\varepsilon)$ if $\varepsilon < 1/6$ ([1], (19.13.1) and (19.9.5)), so they certainly hold for $x = O(\sqrt{n}\log n)$. Therefore, when $s = O(\sqrt{n}\log n)$, the above integral, which equals

$$\int_0^s \frac{1}{2\sqrt{n}} K\left(\frac{x}{2\sqrt{n}}, \frac{x}{2\sqrt{n}}\right) dx,$$

also equals

$$\frac{1}{\pi} \int_0^s \left(1 - \frac{\sin 2x}{2x} + o(1)\right) dx \geq \delta s$$

for some $\delta$. This gives (5.6) and concludes the proof of (5.4).

For the highest curve, we use the top scaling in trace norm of the $m = 1$ kernel, which is just $H_{11}(x,y)$, as described in Section 5.3. We find, for fixed $\tau$, that $\hat{X}_n$ (a random variable independent of $\tau$) defined by

$$\frac{X_n}{\sqrt{2\tau(1-\tau)}} = 2\sqrt{n} + \frac{\hat{X}_n}{2^{2/3} n^{1/6}}$$

has the limit distribution $F_2$, the limiting distribution of the largest eigenvalue in GUE [23]. Since

$$\mathbb{E}(A_{n,H}) = \int_0^1 \mathbb{E}(X_n(\tau))\, d\tau = \frac{\pi}{2^{3/2}} \sqrt{n} + \frac{\pi}{8(2n)^{1/6}} \mathbb{E}(\hat{X}_n),$$

we see that as $n \to \infty$,

(5.7) $$\mathbb{E}(A_{n,H}) = \frac{\pi}{2^{3/2}} \sqrt{n} + \frac{c_H}{n^{1/6}} + o(n^{-1/6}),$$

where

$$c_H = \frac{\pi}{8 \cdot 2^{1/6}} \int_{-\infty}^{\infty} s\, dF_2(s).$$

Here,

$$\int_{-\infty}^{\infty} s\, dF_2(s) = -1.771086807411601\ldots,$$

as computed to high precision by M. Prähofer, so

$$c_H \simeq -0.619623767170\ldots.$$

Evaluating (5.7) for $n = 5, \ldots, 9$ gives 2.00981, 2.26104, 2.49069, 2.70345, 2.90254, respectively. These numbers should be compared with the entries in the third column of Table 1.

**Acknowledgments.** The authors thank Chris Cosgrove for his invaluable assistance regarding the differential equation (1.2) and Neil O'Connell for elucidation of the other Bessel process. The first author thanks Anne Boutet de Monvel for her invitation to visit, and her kind hospitality at, the University of Paris 7, where part of this work was done.

DEPARTMENT OF MATHEMATICS  
UNIVERSITY OF CALIFORNIA  
DAVIS, CALIFORNIA 95616  
USA  
E-MAIL: [tracy@math.ucdavis.edu](mailto:tracy@math.ucdavis.edu)

DEPARTMENT OF MATHEMATICS  
UNIVERSITY OF CALIFORNIA  
SANTA CRUZ, CALIFORNIA 95064  
USA  
E-MAIL: [widom@ucsc.edu](mailto:widom@ucsc.edu)